\documentclass[reqno,11pt]{amsart}
\usepackage{amssymb,amsmath,amsthm,amsfonts}

\usepackage{enumerate,esint}
\usepackage{mathtools}
\usepackage{graphicx}

\usepackage{epstopdf}
\usepackage[dvipsnames]{xcolor}
\usepackage{hyperref}	
\definecolor{ourcolor}{RGB}{0,102,204}

\usepackage[
   rm={oldstyle=false,proportional=true},
   sf={oldstyle=true,proportional=true},
   tt={oldstyle=true,proportional=true,variable=true},
   qt=false,
 ]{cfr-lm}
 
\usepackage[normalem]{ulem}
\usepackage{hyperref}
 \hypersetup{
 colorlinks,
 citecolor=ourcolor,
 linkcolor=ourcolor,
 urlcolor=black}

\numberwithin{equation}{section}

\usepackage{geometry}
 \makeatletter
 \def\@seccntformat#1{\hspace*{0mm}%
  \protect\textup{\protect\@secnumfont
    \ifnum\pdfstrcmp{subsection}{#1}=0 \bfseries\fi
    \csname the#1\endcsname
    \protect\@secnumpunct
      }%
 }

\usepackage[normalem]{ulem}

\newcommand{\assign}{:=}
\newcommand{\mathd}{\mathrm{d}}
\newcommand{\of}{:}
\newcommand{\tmSep}{; }
\newcommand{\tmcolor}[2]{{\color{#1}{#2}}}
\newcommand{\tmem}[1]{{\em #1\/}}
\newcommand{\tmmathbf}[1]{\ensuremath{\boldsymbol{#1}}}
\newcommand{\tmname}[1]{\textsc{#1}}

\newcommand{\tmsep}{, }
\newcommand{\tmtextsf}[1]{\text{{\sffamily{#1}}}}

{\theoremstyle{remark}\newtheorem{remark}{Remark}}
\newtheorem{theorem}{Theorem}

\newcommand{\Aop}{K}
\newcommand{\aeps}{\alpha_{\varepsilon}}

\newcommand{\zv}{z}

\newcommand{\vv}{\tmmathbf{v}}

\newcommand{\jM}{\mathbb{S}^2}

\newcommand{\Om}{\text{\tmtextsf{\ensuremath{\Omega}}}}
\newcommand{\T}{\scriptstyle{\top}}

\newcommand{\RR}{\mathbb{R}}
\newcommand{\NN}{\mathbb{N}}
\newcommand{\eqs}{=}

\begin{document}

 \title[Exchange interactions in curved magnetic nanowires]{Reduced theory of symmetric and antisymmetric exchange interactions in nanowires}

\author{G{\tmname{iovanni}} D{\tmname{i}} {\tmname{Fratta}}}
\address{Giovanni Di Fratta, Dipartimento di Matematica e Applicazioni ``R.
Caccioppoli'', Universit{\`a} degli Studi di Napoli ``Federico II'', Via
Cintia, Complesso Monte S. Angelo, 80126 Napoli, Italy}

\author{{\tmname{Filipp N. Rybakov}}}
\address{Filipp N. Rybakov, Department of Physics and Astronomy, Uppsala University, Box-516, Uppsala SE-751 20, Sweden}

\author{{\tmname{Valeriy Slastikov}}}
\address{Valeriy Slastikov, School of Mathematics, University of\\
Bristol, University Walk, Bristol BS8 1TW, United Kingdom}

\begin{abstract}
   We investigate the behavior of minimizers of perturbed Dirichlet energies supported on a wire generated by a regular simple curve $\gamma$ and defined in the space of $\mathbb{S}^2$-valued functions. The perturbation $K$ is represented by a matrix-valued function defined on $\mathbb{S}^2$ with values in $\RR^{3 \times 3}$. Under natural regularity conditions on $K$, we show that the family of perturbed Dirichlet energies converges, in the sense of $\Gamma$-convergence, to a simplified energy functional on $\gamma$. The reduced energy unveils how part of the antisymmetric exchange interactions contribute to an anisotropic term whose specific shape depends on the curvature of $\gamma$.
   
We also discuss the significant implications of our results for studies of ferromagnetic nanowires when Dzyaloshinskii-Moriya interaction (DMI) is present. 
\end{abstract}



\maketitle

\section{Introduction and motivation}

The motivation for our study comes from investigating the behavior of ferromagnetic nanowires.
Microscopically ferromagnetic nanowires such as nanotubes, nanowhiskers, or nanostrips can exhibit unique physical properties. These include anisotropic magnetoresistance in helimagnets~{\cite{Du2014,Mathur2019}}, proximity exchange fields in ferromagnet--superconductor heterostructures~{\cite{Liu2020,Vaitiekenas2021}}, fast domain wall motion~{\cite{Schoebitz2019}}, impact of wire twisting onto magnetic textures~{\cite{SanzHernandez2020}}, and many
others~{\cite{FernandezPacheco2017,Makarov2022}}. 
Techniques for fabricating magnetic nanowires make available a wide variety of cross-sections and shapes, from almost perfectly straight to
braided~{\cite{Huang2021,SanzHernandez2020}}. Of significant interest for research and applications are wireframes and networks made of magnetic nanowires~{\cite{Sahoo2018,Keller2018,Gliga2019,Volkov2024}}.

The theoretical description of the statics and dynamics of magnetization in
nanomagnets is a challenging task even from the numerical point of view due to the nonlinear, nonconvex, and nonlocal nature of the problem, see, e.g., the monograph~{\cite{prohl2001}}, the review
articles~{\cite{KruzikProhl06,carlos2007}} and the references therein. Brute-force numerical simulations are typically performed within the framework of 3D micromagnetic model~{\cite{Exl2019,Zheng2021,Cheenikundil2021,Volkov2024}}, which are
challenging and often time-consuming.

While an essential goal of quantitative simulations is description,
the basic explanation of magnetic phenomena typically requires having simplified yet accurate models at
hand.
From the mathematical point of view, such models can be derived by dimension reduction techniques of the calculus of variations as $\Gamma$-limits of the micromagnetic energy functional under suitable asymptotic relations among the material and geometric scaling parameters.
Dimensional reduction arguments in micromagnetics are a subject with a long history as they date back to the seminal paper {\cite{GioiaJames97}}, where the authors show that in
{\emph{planar}} thin films, the nonlocal effects of the stray field operator
reduce to an easy-plane anisotropy term. Later, various static reduced theories for planar thin-film micromagnetics were established under different scaling regimes~ {\cite{carbou2001thin,DeSimoneKohnMuellerOtto02,DeSimoneETal01,KohnSlastikov05,Moser04,Moser05,Slastikov05}}.
Further asymptotic regimes, e.g., in the presence of strong perpendicular
anisotropy, have been the subject of~{\cite{Di_Fratta_2024}}.
The induced thin-film dynamics have been analyzed in~{\cite{KohnSlastikov-dyn05,Melcher10,EGarcaCervera2001EffectiveDF,Capella_2007}}.

The latest developments in nanotechnology have made it possible to create nanostructures, with thicknesses as small as a few atomic layers and lateral sizes down to tens of nanometers.
These structures often demonstrate the impact of interfacial effects, such as the {\emph{Dzyaloshinskii-Moriya interaction}} (DMI)~{\cite{Dzyaloshinsky_1958,Moriya_1960}}, leading to the formation of magnetic skyrmions {\cite{fert2017magnetic}}. Mathematical studies of magnetic domain walls in the regime of planar thin films, taking into account also DMI, have been done in {\cite{Davoli_2022,melcher14,ms:prsla17,mskt:prb17}}.

In recent years, interest has grown in ferromagnetic systems with a {\emph{curved}}  shape due to their capability of hosting magnetic skyrmions even in materials where the DMI can be neglected. The evidence of these states sheds light on the role of the geometry in magnetism: chiral spin-textures can be stabilized by curvature effects only, in contrast to the planar case where DMI is required (we refer the reader to the recent monograph {\cite{Makarov2022}}). From the mathematical perspective, dimension reduction results for curved thin films have been
studied in~\cite{carbou2001thin,Di_Fratta_2020,doi:10.1137/19M1261365,Slastikov05}
and, more recently, in {\cite{Di_Fratta_2023}} where the curved thin film
limit of chiral Dirichlet energies is derived.

The effective magnetic behavior of ferromagnetic nanowires has been investigated in {\cite{Slastikov2012}} where both long-range dipolar interactions and symmetric exchange interactions are considered, also for curved wire. The main aim of this paper is to complement the analysis in {\cite{Slastikov2012}} by taking into account general antisymmetric exchange interactions. As in {\cite{Davoli_2022}} and {\cite{Di_Fratta_2023}}, the derived limiting model reveals new physics: part of the antisymmetric exchange interactions contributes to an increase in the shape anisotropy originating from the magnetostatic self-energy. Several examples of applications to the micromagnetic theory complete our analytical derivation.

\section{Contributions of the present work}\label{sec:sec2}
\begin{figure}[t]
	\centering
	\includegraphics[width=15cm]{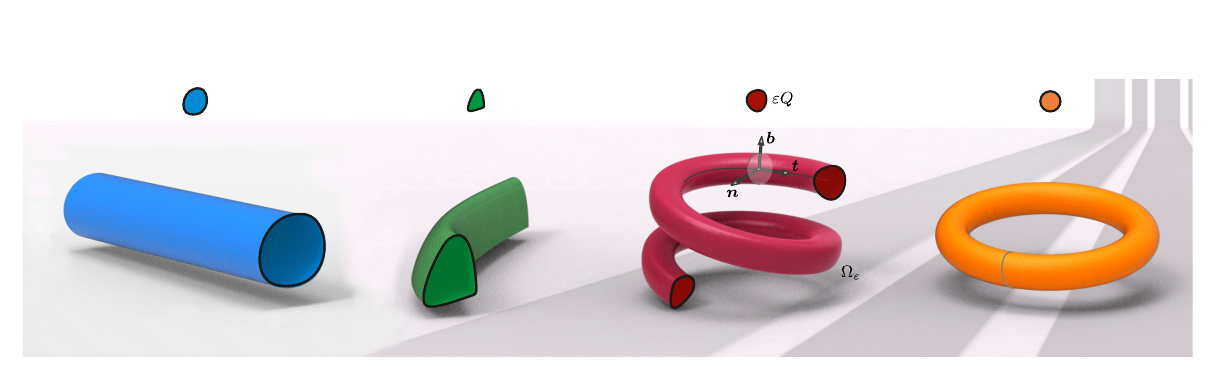}
	\caption{\small Sketches of suitable wires. The wires can be straight or curved and can have a nonsmooth boundary or a multiply connected cross-section (not every case is shown in the picture). The cross-section is assumed to be invariant along the wire.}
	\label{fig_wires}
\end{figure}

Given a regular simple curve\footnote{By a regular simple curve, we mean the
image of a $C^2$-map $\gamma \of I \mapsto \RR^3$, $I \subseteq \RR$ a compact
interval, such that $\partial_s \gamma (s) \neq 0$ for every $s \in I$, and
with no self-intersections, i.e., such that the only possible loss of
injectivity in $\gamma$ arises at the endpoints of $I$, case in which the
curve closes into a loop.} $\gamma \of I \mapsto \RR^3$ and a smooth domain $Q
\subseteq \RR^2$, for every $\varepsilon > 0$, we denote by $\Om_{\varepsilon}
\subseteq \RR^3$ the $\varepsilon$-tube along $\gamma$ with cross-sectional
shape $Q$ and thickness $\varepsilon$, defined as the image of the cylindrical
region $I \times Q \subseteq \RR \times \RR^2$ through the parameterization
(see Figure \ref{fig_wires})
\begin{equation}
  \varphi_{\varepsilon} \of \left( s, \zv \right) \in I \times Q \mapsto
  \gamma (s) + \varepsilon z_1 \tmmathbf{n} (s) + \varepsilon z_2 \tmmathbf{b}
  (s) \in \Om_{\varepsilon} . \label{eq:epstube}
\end{equation}
Here, we set $\zv \assign (z_1, z_2)$ and denoted by $(\tmmathbf{t}(s),
\tmmathbf{n}(s), \tmmathbf{b}(s))$ the {\tmname{Frenet--Serret}} frame
consisting, respectively, of the tangent, normal, and binormal unit vectors to
the curve at $\gamma (s)$ --- if $\gamma$ is a straight line, we complete
\ensuremath{\tmmathbf{t}} to an orthonormal basis of $\RR^3$. In what follows,
we assume that $\gamma$ is parameterized by arc length. Also, we assume that
the interval $I$ is compact so that according to the tubular neighborhood theorem, there exists a $\delta > 0$ such that for every $0 < \varepsilon < \delta$, the map $\varphi_\varepsilon$ defined by \eqref{eq:epstube} is a diffeomorphism of $I\times Q$ onto the $\varepsilon$-tube along $\gamma$.

The main aim of this paper is to investigate the curved thin-wire limit
($\varepsilon \rightarrow 0$) of the perturbed Dirichlet energy ($0 <
\varepsilon < \delta$)
\begin{equation}
  \mathcal{G}_{\varepsilon}( \vv ) \assign \frac{1}{2
  \varepsilon^2 | Q |} \int_{\Om_{\varepsilon}} \left| D \vv (x) + \Aop (
  \vv (x) ) \right|^2 \mathd x \label{eq:mainenfunc}
\end{equation}
defined on $H^1$ Sobolev maps $\vv : \Om_{\varepsilon} \rightarrow
\mathbb{S}^2$, where for every $\varepsilon > 0$ the domain $\Om_{\varepsilon}
\subseteq \RR^3$ is an $\varepsilon$-tube of $\RR^3$, and $\mathbb{S}^2$ is
the two-sphere of $\RR^3$. In writing \eqref{eq:mainenfunc}, we made the
common and convenient abuse of denoting by $| \cdot |$ both the Euclidean norm
on $\RR^{3 \times 3}$ and the Lebesgue measure when applied to sets as in the
expression $| Q |$.

The perturbation $\Aop : \sigma \in \mathbb{S}^2 \mapsto \Aop (\sigma) \in
\mathbb{R}^{3 \times 3}$ is represented by a matrix-valued function defined on
$\mathbb{S}^2$ and with values in $\RR^{3 \times 3}$ and is assumed to be
Lipschitz continuous, i.e., that there exists $c_{\Aop} > 0$ such that
\begin{equation}
  \left| \Aop (\sigma_1) - \Aop (\sigma_2) \right|_{n \times m} \leqslant
  c_{\Aop} | \sigma_1 - \sigma_2 |_m \quad \forall \sigma_1, \sigma_2 \in
  \mathbb{S}^2 \label{eq:cALip} .
\end{equation}
Since $\mathbb{S}^2$ is compact, the continuity condition \eqref{eq:cALip}
implies that the image of $\Aop$ is bounded (in fact, compact). In what
follows, to simplify notation, we will assume that a Lipschitz constant $c_{\Aop}$
is chosen big enough so that there holds
\begin{equation}
  \left| \Aop (\sigma) \right| \leqslant c_{\Aop} \quad \forall \sigma \in
  \mathbb{S}^2 \label{eq:cALipnew} .
\end{equation}

\noindent For any $0 < \varepsilon < \delta$, the existence of at least a minimizer for
$\mathcal{G}_{\varepsilon}$ in $H^1( \Om_{\varepsilon}, \mathbb{S}^2)$ is a simple application of the direct method of the calculus of
variations. We are interested in the asymptotic behavior of the family of
minimizers of $(\mathcal{G}_{\varepsilon})_{\varepsilon}$ as $\varepsilon
\rightarrow 0$. Our main result ({cf.} Theorem~\ref{thm:main} below)
shows that the family $(\mathcal{G}_{\varepsilon})_{\varepsilon}$ converges,
in the sense of $\Gamma$-convergence to an energy functional defined on the
curve $\gamma$, which strongly depends on $K$ and the target space
$\mathbb{S}^2$, and has remarkable physical implications for relevant systems,
for instance, to magnetic materials: in the curved thin-wire regime, generic
exchange interactions (symmetric and antisymmetric) can manifest themselves
under an additional anisotropy term whose specific shape depends both on the
shape of the wire and the curvature of the target surface.

To properly state our result, we must finish setting up the stage. For every
$0 < \varepsilon < \delta$\quad we introduce the following functional defined
on $H^1 (\mathcal{C}, \mathbb{S}^2)$, which can be thought of as the pull-back
of $\mathcal{G}_{\varepsilon}$ on the cylindrical domain $\mathcal{C} \assign
I \times Q$:
\begin{align}
  \mathcal{E}_{\mathcal{C}}^{\varepsilon} \left( \vv_{\varepsilon} \right) &
  \assign  \frac{1}{2 | Q |} \int_{I \times Q} \left| \partial_s
  \vv_{\varepsilon} + \tau \left( \zv \wedge \nabla_{\zv} \vv_{\varepsilon}
  \right) + \aeps \Aop ( \vv_{\varepsilon}) \tmmathbf{t} \right|^2
  \frac{1}{\aeps} \mathd \zv \mathd s \nonumber\\
  &   \; + \frac{1}{2 | Q |} \int_{I \times Q} \left( |
  \varepsilon^{- 1} \partial_1 \vv_{\varepsilon} + \aeps \Aop \left(
  \vv_{\varepsilon} \right) \tmmathbf{n}|^2 + | \varepsilon^{- 1}
  \partial_2 \vv_{\varepsilon} + \aeps \Aop ( \vv_{\varepsilon} )
  \tmmathbf{b} |^2 \right)  \frac{1}{\aeps} \mathd \zv \mathd s 
  \label{eq:rewriteG1}
\end{align}
where we set
\begin{equation} \label{eqs:defsnot}
\partial_1 \assign \partial_{z_1}, \partial_2 \assign \partial_{z_2},
   \quad \zv \wedge \nabla_{\zv} \vv_{\varepsilon} \assign z_2 \partial_1
   \vv_{\varepsilon} - z_1 \partial_2 \vv_{\varepsilon}, \quad \aeps \left(
   s, \zv \right) \assign (1 - \varepsilon \varpi(s) z_1),
\end{equation}
with $\varpi(s)$ being the curvature of the space curve $\gamma$ at $\gamma (s)$.
Our main result is stated in the following statement.
\begin{theorem}
  \label{thm:main}For any $0<\varepsilon<\delta$, the minimization
  problem for $\mathcal{G}_{\varepsilon}$ in $H^1 (\Om_{\varepsilon},
  \mathbb{S}^2)$ is equivalent to the minimization in $H^1
  (\mathcal{C}, \mathbb{S}^2)$ of the functional
  $\mathcal{E}_{\mathcal{C}}^{\varepsilon}$ defined by
  {\emph{\eqref{eq:rewriteG1}}} in the sense that a configuration
  $\vv_{\varepsilon} \in H^1 ( \Om_{\varepsilon}, \jM )$ minimizes
  $\mathcal{G}_{\varepsilon}$ if and only if $\vv_{\varepsilon} \circ
  \varphi_{\varepsilon} \in H^1 (\mathcal{C}, \mathbb{S}^2)$ minimizes
  $\mathcal{E}_{\mathcal{C}}^{\varepsilon}$.
  
  The family $(\mathcal{E}_{\mathcal{C}}^{\varepsilon})_{\varepsilon }$ is equicoercive in the weak topology of $H^1 (\mathcal{C},
  \mathbb{S}^2)$ and the $\Gamma$-limit $\mathcal{E}_{\mathcal{C}} \assign
  \Gamma \text{-} \lim_{\varepsilon \rightarrow 0}
  \mathcal{E}_{\mathcal{C}}^{\varepsilon}$ is defined for every $\tmmathbf{v}
  \in H^1( \mathcal{C}, \jM)$ by
  \begin{equation}
    \mathcal{E}_{\mathcal{C}}( \vv ) = \left\{ \begin{array}{ll}
      \displaystyle \frac{1}{2} \int_I | \partial_s \vv (s) + \Aop( \vv (s)) \tmmathbf{t} (s) |^2 \mathd s & \\
      \qquad + \displaystyle \frac{1}{2} \int_I | \Aop^{\T} (\vv (s)) \vv
      (s) |^2 - ( \Aop^{\T} (\vv (s)) \vv (s) \cdot
      \tmmathbf{t} (s) )^2 \mathd s & \text{\quad if } {\nabla_{\zv}}  \vv =
      0,\\
      + \infty & \text{\quad otherwise} .
    \end{array} \right. \label{eq:gleF}
  \end{equation}
  Moreover,
  \begin{equation}
    \min_{H^1 \left( \Om_{\varepsilon}, \jM \right)} \mathcal{G}_{\varepsilon}
    = \min_{H^1( \mathcal{C}, \jM)}
    \mathcal{E}_{\mathcal{C}}^{\varepsilon} {= \min_{H^1 \left( \mathcal{C},
    \jM \right)}}  \mathcal{E}_{\mathcal{C}} + o (1)
    \label{eq:firstordergammadevelopforF}
  \end{equation}
  and if $( \vv_{\varepsilon})_{\varepsilon}$ is
  a minimizing family for
$(\mathcal{E}_{\mathcal{C}}^{\varepsilon})_{\varepsilon}$
  then $(\vv_{\varepsilon})_{\varepsilon}$
  converges, strongly in $H^1( \mathcal{C}, \jM)$, to a minimum
  point of $\mathcal{E}_{\mathcal{C}}$.
\end{theorem}

\begin{remark} \label{rmk:altdescgammlim}
  In terms of normal and binormal vectors to the curve, the expression in
  \eqref{eq:gleF} is equivalent to
  \[ \frac{1}{2} \int_I \left| \partial_s \vv (s) + \Aop(\vv (s)) \tmmathbf{t} (s) \right|^2 \mathd s + \frac{1}{2} \int_I (\Aop^{\T}(\vv (s)) \vv (s) \cdot \tmmathbf{n} (s))^2
     + ( \Aop^{\T}(\vv (s)) \vv (s) \cdot \tmmathbf{b} (s))^2 \mathd s. \]
\end{remark}

\begin{proof}
\tmname{(Equivalence)} First, we show that for any $0<\varepsilon<\delta$, the minimization
  problem for $\mathcal{G}_{\varepsilon}$ in $H^1 \left( \Om_{\varepsilon},
  \mathbb{S}^2 \right)$ is equivalent to the minimization in $H^1
  (\mathcal{C}, \mathbb{S}^2)$ of the functional
  $\mathcal{E}_{\mathcal{C}}^{\varepsilon}$ defined by
  {\emph{\eqref{eq:rewriteG1}}}. For that, let us denote by with $F \assign (\tmmathbf{t} |
  \tmmathbf{n}|\tmmathbf{b})$ the matrix whose columns are the component of
  the Frenet frame. By the Frenet--Serret formulas (with $\varpi$ and $\tau$ denoting, respectively, the curvature and the torsion of $\gamma$) :
  \[ \dot{\tmmathbf{n}} (s) = - \varpi \tmmathbf{t}+ \tau \tmmathbf{b}, \quad
     \dot{\tmmathbf{t}} = \varpi \tmmathbf{n}, \quad \dot{\tmmathbf{b}} = -
     \tau \tmmathbf{n}, \]
  we get that for every $0 < \varepsilon < \delta$, there holds
  \[ D \varphi_{\varepsilon} = \left(\begin{array}{c|c|c}
       \partial_s \varphi_{\varepsilon} & \partial_1 \varphi_{\varepsilon} &
       \partial_2 \varphi_{\varepsilon}
     \end{array}\right) = F \cdot \left(\begin{array}{ccc}
       1 - \varepsilon \varpi z_1 & 0 & 0\\
       - \varepsilon \tau z_2 & \varepsilon & 0\\
       \varepsilon \tau z_1 & 0 & \varepsilon
     \end{array}\right) \]
  and $\det D \varphi_{\varepsilon} = \varepsilon^2 \aeps$, with $\aeps (s, z_1)
     \assign 1 - \varepsilon \varpi (s) z_1$.
  Note that $D \varphi_{\varepsilon}$ is invertible for $\varepsilon$ small
  because. The inverse of $D \varphi_{\varepsilon}$ is given by
  \[ (D \varphi_{\varepsilon})^{- 1} = \frac{1}{\aeps} \Phi_{\varepsilon}
     F^{\T}, \quad \Phi_{\varepsilon} \assign \left(\begin{array}{ccc}
       1 & 0 & 0\\
       \tau z_2 & \varepsilon^{- 1} & 0\\
       - \tau z_1 & 0 & \varepsilon^{- 1}
     \end{array}\right) . \]
  The previous considerations allow showing the equivalence between the
  minimization problem for $\mathcal{G}_{\varepsilon}$ and
  $\mathcal{E}_{\mathcal{C}}^{\varepsilon}$. Indeed with $\vv_{\varepsilon}
  \assign \vv \circ \varphi_{\varepsilon}$, we have
  \begin{align}
    \mathcal{G}_{\varepsilon} \left( \vv \right) &\eqs  \frac{1}{2
    \varepsilon^2 | Q |} \int_{\Om_{\varepsilon}} \left| D \vv (x) + \Aop
    \left( \vv (x) \right) \right|^2 \mathd x \nonumber\\
    & =  \frac{1}{2} \int_{I \times Q} \left| D \vv_{\varepsilon}  (D
    \varphi_{\varepsilon})^{- 1} + \Aop (\vv_{\varepsilon})
    \right|^2  \aeps \mathd \zv \mathd s \nonumber\\
   & =  \frac{1}{2 | Q |} \int_{I \times Q} \left| \left( \partial_s
    \vv_{\varepsilon} + \tau \left( z_2 \partial_1 \vv_{\varepsilon} - z_1
    \partial_2 \vv_{\varepsilon} \right) \right) + \aeps \Aop \left(
    \vv_{\varepsilon} \right) \tmmathbf{t} \right|^2  \frac{1}{\aeps} \mathd
    \zv \mathd s \nonumber\\
    & \qquad \qquad + \frac{1}{2 | Q |} \int_{I \times Q} \left( | \varepsilon^{- 1}
    \partial_1 \vv_{\varepsilon} + \aeps \Aop( \vv_{\varepsilon} )
    \tmmathbf{n} |^2 + | \varepsilon^{- 1} \partial_2
    \vv_{\varepsilon} + \aeps \Aop ( \vv_{\varepsilon})
    \tmmathbf{b} |^2 \right)  \frac{1}{\aeps} \mathd \zv \mathd s
    \nonumber\\
    &=  \mathcal{E}_{\mathcal{C}}^{\varepsilon} (\tmmathbf{v}_{\varepsilon}),
    \nonumber
  \end{align}
  from which the equivalence follows.
  
   \smallskip
   
  {\noindent}{\tmname{(Compactness/Equicoerciveness)}} 
We want to show that the family
$(\mathcal{E}_{\mathcal{C}}^{\varepsilon})_{\varepsilon}$ is equicoercive in
the weak topology of $H^1 (\mathcal{C}, \mathbb{S}^2)$, i.e., the existence of
a weakly compact set $J (\mathcal{C}, \mathbb{S}^2) \subseteq H^1
(\mathcal{C}, \mathbb{S}^2)$ such that for every $0 < \varepsilon < \delta$
there holds
\begin{equation}
  \min_{H^1 (\mathcal{C}, \mathbb{S}^2)}
  \mathcal{E}_{\mathcal{C}}^{\varepsilon} = \min_{J (\mathcal{C},
  \mathbb{S}^2)} \mathcal{E}_{\mathcal{C}}^{\varepsilon} .
  \tmcolor{red}{\label{eq:equicoerciveE}}
\end{equation}
Ensuring the equicoercivity of
$(\mathcal{E}_{\mathcal{C}}^{\varepsilon})_{\varepsilon}$ will guarantee that
the fundamental theorem of $\Gamma$-convergence hypotheses are met.
Consequently, the $\Gamma$-limit $\mathcal{E}_{\mathcal{C}}$ will satisfy
realtion \eqref{eq:firstordergammadevelopforF} on the variational convergence
of minimum problems.

To show \eqref{eq:equicoerciveE}, we first observe that for any constant in
space configuration $\xi \in H^1 (\mathcal{C}, \mathbb{S}^2)$ we have
\begin{equation}
  \min_{\vv \in H^1 (\mathcal{C}, \mathbb{S}^2)}
  \mathcal{E}_{\mathcal{C}}^{\varepsilon} \left( \vv \right) \leqslant
  \mathcal{E}_{\mathcal{C}}^{\varepsilon} (\xi) = \frac{1}{2 |Q|}  \int_{I
  \times Q}  \frac{\left| \Aop (\xi) \right|^2}{\sqrt{\aeps ( s, \zv
  )}} \mathd \zv \mathd s.
\end{equation}
Given that for $\varepsilon$ sufficiently small $\aeps \left( s, \zv \right)$
is uniformly bounded away from zero, there exists $c_{\varpi} > 0$ such that
(cf. the expression of $\alpha_{\varepsilon}$ in \eqref{eqs:defsnot})
\begin{equation}
  c_{\varpi}^{- 1} < \aeps \left( s, \zv \right) < c_{\varpi} \quad \text{ for
  every } \left( s, \zv \right) \in I \times Q.
\end{equation}
Also, as stated in \eqref{eq:cALipnew}, $\left| \Aop (\sigma) \right|
\leqslant c_{\Aop}$. Therefore, for $\delta$ sufficiently small, the
minimizers of $(\mathcal{E}_{\mathcal{C}}^{\varepsilon})_{0 < \varepsilon <
\delta}$ are all contained in the $\varepsilon$-independent set
\begin{equation}
  J (\mathcal{C}, \mathbb{S}^2) \assign \bigcup_{0 < \varepsilon < \delta}
  \left\{ \vv \in H^1 (\mathcal{C}, \mathbb{S}^2) \of
  \mathcal{E}_{\mathcal{C}}^{\varepsilon} \left( \vv \right) \leqslant c
  \right\}, \label{eq:setJ}
\end{equation}
where $c > 0$ is a positive constant that depends only on $c_{\varpi},
c_{\Aop}$ and on the geometry of $\gamma$. After that, we use that for every
$a, b \in \mathbb{R}$ and every $\eta > 0$ there holds $|a + b|^2 \geqslant (1
+ \eta)^{- 1} |a|^2 - \eta^{- 1} |b|^2$, to infer that
\begin{eqnarray}
  2 |Q| \mathcal{E}_{\mathcal{C}}^{\varepsilon} \left( \vv \right) & \geqslant
  & \frac{1}{c_{\varpi} (1 + \eta)}  \int_{I \times Q} \left| \partial_s  \vv
  + \tau \left( \zv \wedge \nabla_{\zv}  \vv \right) \right|^2 \nonumber\\
  &  & \quad + \frac{1}{\aeps} \mathd \zv \mathd s + \frac{1}{\varepsilon^2
  (1 + \eta)} \int_{I \times Q} | \nabla_{\zv}  \vv |^2  \frac{1}{\aeps}
  \mathd \zv \mathd s - \frac{1}{\eta} \int_{I \times Q}  \frac{\left| \Aop
  \left( \vv \right)  \right|^2}{\sqrt{\aeps}} \mathd \zv \mathd s \nonumber\\
  & \geqslant & \frac{1}{c_{\varpi} (1 + \eta)^2}  \int_{I \times Q} \left|
  \partial_s  \vv \right|^2 \mathd \zv \mathd s - \frac{1}{c_{\varpi} \eta (1
  + \eta)} \tau^2 \int_{I \times Q} \left|  \zv \wedge \nabla_{\zv}  \vv
  \right|^2 \mathd \zv \mathd s \nonumber\\
  &  & \quad + \frac{1}{c_{\varpi} (1 + \eta) \varepsilon^2} \int_{I \times
  Q} | \nabla_{\zv}  \vv |^2  \frac{1}{\aeps} \mathd \zv \mathd s -
  \frac{c_{\Aop}^2}{\eta \sqrt{c_{\varpi}}} | I \times Q | . \nonumber
\end{eqnarray}
Given that $\left| \zv \wedge \nabla_{\zv} \vv \right| \leqslant 2 | z | \cdot
\left| \nabla_{\zv} \vv \right| \leqslant c_Q \left| \nabla_{\zv} \vv \right|$
for a positive constant $c_Q > 0$ that depends only on $Q$, we get that up to
the constant term $- \frac{c_{\Aop}^2}{\eta \sqrt{c_{\varpi}}} | I \times Q |$
there holds
\[ 2 |Q| \mathcal{E}_{\mathcal{C}}^{\varepsilon} \left( \vv \right) \geqslant
   \frac{1}{c_{\varpi} (1 + \eta)^2} \int_{I \times Q} \left| \partial_s  \vv
   \right|^2 \mathd \zv \mathd s + \frac{1}{(1 + \eta) c_{\varpi} } \left(
   \frac{1}{\varepsilon^2} - \frac{\tau^2 c_Q^2}{\eta} \right) \int_{I \times
   Q} | \nabla_{\zv} \vv |^2 \mathd \zv \mathd s. \]
In particular, choosing $\eta \assign 2 \varepsilon^2 \tau^2 c_Q^2$, we get
the existence of a constant $c_{\ast} > 0$ dependent only on $\delta,
c_{\varpi}, c_{\Aop}$ and the shape of the $\delta$-tube generated by
$\gamma$, such that
\begin{equation}
  \int_{I \times Q} \left| \partial_s \vv \right|^2 \mathd \zv \mathd s +
  \frac{1}{\varepsilon^2} \int_{I \times Q} | \nabla_{\zv}  \vv_{\varepsilon}
  |^2 \mathd \zv \mathd s \leqslant c_{\ast} \left( 1
  +\mathcal{E}_{\mathcal{C}}^{\varepsilon} \left( \vv \right) \right) .
  \label{eq:est4compact}
\end{equation}
Assuming without loss of generality that $\delta < 1$ so that also
$\varepsilon < 1$, from the previous relation \eqref{eq:est4compact} we can
conclude that the set $J (\mathcal{C}, \mathbb{S}^2)$ defined by
\eqref{eq:setJ} is contained in the bounded subset of $H^1 (\mathcal{C},
\mathbb{R}^3)$ given by
\[ H^1_b (\mathcal{C}, \mathbb{S}^2) \assign \{ \tmmathbf{v} \in H^1
   (\mathcal{C}, \mathbb{R}^3) \of | \tmmathbf{v} | = 1, \| \tmmathbf{v}
   \|^2_{H^1 (\mathcal{C}, \mathbb{R}^3)} \leqslant c_{\ast} (1 + c) \} \]
Therefore
\[ \min_{\vv \in H^1 (\mathcal{C}, \mathbb{S}^2)}
   \mathcal{E}_{\mathcal{C}}^{\varepsilon} \left( \vv \right) = \min_{\vv \in
   H^1_b (\mathcal{C}, \mathbb{S}^2)} \mathcal{E}_{\mathcal{C}}^{\varepsilon}
   \left( \vv \right) . \]
To establish that $H^1_b (\mathcal{C}, \mathbb{S}^2)$ is weakly compact, we
need to prove that it is weakly closed. To do this, consider a sequence \
$\left( \vv_n \right)_{n \in \mathbb{N}}$ in $H^1_b (\mathcal{C},
\mathbb{S}^2)$ such that $\vv_n \rightharpoonup \vv_0$ weakly in $H^1
(\mathcal{C}, \mathbb{R}^3)$; by Rellich-Kondrachov theorem, $\vv_n
\rightarrow \vv_0$ strongly in $L^2 (\mathcal{C}, \mathbb{R}^3)$. Thus, after
extracting a subsequence that converges pointwise a.e. in $\mathcal{C}$, we
get that the $\vv_0$ still takes values in $\mathbb{S}^2$. This completes the
proof of equicoerciveness.

      \smallskip
      
  {\noindent}{\tmname{($\Gamma$-liminf)}} Let 
  $(\tmmathbf{v}_{\varepsilon})$ be a family in $H^1 (\mathcal{C},
  \mathbb{S}^2)$ such that
  \[ \sup_{0 < \varepsilon < \delta} \mathcal{E}_{\mathcal{C}}^{\varepsilon}
     (\tmmathbf{v}_{\varepsilon}) < + \infty . \]
  By estimate \eqref{eq:est4compact}, we get that the families $\left( \varepsilon^{- 1}
  \partial_1 \vv_{\varepsilon} \right)$, $\left( \varepsilon^{- 1} \partial_2
  \vv_{\varepsilon} \right)$ and $\left( \partial_s \vv_{\varepsilon} \right)$
  are bounded in $L^2$. Therefore, there exist $\tmmathbf{v}_0 \in H^1
  (\mathcal{C}, \mathbb{S}^2), \tmmathbf{d}_1, \tmmathbf{d}_2 \in L^2 (
  \mathcal{C}, \RR^3)$, such that
  \begin{equation} \label{eqs:compact1}
  \arraycolsep=1.4pt\def\arraystretch{1.2}
  \begin{array}{rcll}
    \tmmathbf{v}_{\varepsilon} & \rightarrow & \tmmathbf{v}_0  
    & \text{\quad strongly in } L^2 (\mathcal{C}, \mathbb{S}^2), \\
    \partial_s \tmmathbf{v}_{\varepsilon} & \rightharpoonup & \partial_s
    \tmmathbf{v}_0  & \text{\quad weakly in } L^2( \mathcal{C}, \RR^3),
   \\
    \frac{1}{\varepsilon} \partial_1 \vv_{\varepsilon} & \rightharpoonup &
    \tmmathbf{d}_1  & \text{\quad weakly in } L^2( \mathcal{C}, \RR^3),
   \\
    \frac{1}{\varepsilon} \partial_2 \vv_{\varepsilon} & \rightharpoonup &
    \tmmathbf{d}_2  & \text{\quad weakly in } L^2( \mathcal{C}, \RR^3) .
  \end{array}
 \end{equation}
  The previous relations imply that $\tmmathbf{v}_0$ is 0-homogeneous along
  the cross-section of the wire, i.e., that it depends only on the
  $s$-variable. In particular,
  \begin{equation} \label{eqs:compact2}
   \nabla_{\zv} \vv_{\varepsilon} \rightarrow 0 \quad \text{strongly in } L^2( \mathcal{C}, \RR^3) .
   \end{equation}
  Moreover, since $\partial_1 \vv_{\varepsilon} \cdot
  \tmmathbf{v}_{\varepsilon} = \partial_2 \vv_{\varepsilon} \cdot
  \tmmathbf{v}_{\varepsilon} = 0$ for every $0 < \varepsilon < \delta$, we
  also infer that $\tmmathbf{d}_1, \tmmathbf{d}_2$ are pointwise orthogonal to
  $\tmmathbf{v}_0$, i.e.,
  \begin{equation} \label{eq:ortconds}
   \tmmathbf{d}_1 \cdot \tmmathbf{v}_0 = 0, \quad \quad \tmmathbf{d}_2 \cdot
     \tmmathbf{v}_0 = 0.
     \end{equation}
  By \eqref{eqs:compact1}, \eqref{eqs:compact2}, and \eqref{eq:ortconds}, taking into account the lower semicontinuity of the norm and that
  $\aeps \rightarrow 1$ when $\varepsilon \rightarrow 0$, we get that
  \begin{align}
    \liminf_{\varepsilon \rightarrow 0}
    \mathcal{E}_{\mathcal{C}}^{\varepsilon} (\tmmathbf{v}_{\varepsilon}) &
    \geqslant  \frac{1}{2} \int_I \left| \partial_s \vv_0 + \Aop \left( \vv_0
    \right) \tmmathbf{t} \right|^2 \mathd s \nonumber\\
    &   + \frac{1}{2 | Q |} \int_{I \times Q} \left(_{_{_{}}} \left|
    \tmmathbf{d}_1 + \Aop (\vv_0) \tmmathbf{n} \right|^2 + \left|
    \tmmathbf{d}_2 + \Aop (\vv_0) \tmmathbf{b} \right|^2 \right)
    \mathd \zv \mathd s \nonumber\\
    & \geqslant  \frac{1}{2} \int_I \left| \partial_s \vv_0 + \Aop \left(
    \vv_0 \right) \tmmathbf{t} \right|^2 \mathd s + \frac{1}{2} \int_I \left[
    \min_{\tmmathbf{d}_1 (s), \tmmathbf{d}_2 (s) \in \vv_0^{\bot} (s)}
    \tmmathbf{g} (\tmmathbf{d}_1 (s), \tmmathbf{d}_2 (s)) \right] \mathd \zv
    \mathd s, \nonumber
  \end{align}
  with 
   \begin{equation} \label{eq:exprg}
  \tmmathbf{g} (\tmmathbf{d}_1, \tmmathbf{d}_2) \assign \left|
     \tmmathbf{d}_1 +  \Aop (\vv_0)
     \tmmathbf{n} \right|^2 + \left| \tmmathbf{d}_2 + 
     \Aop (\vv_0) \tmmathbf{b} \right|^2\, .
     \end{equation}
  The minimal energy of $\tmmathbf{g}$ is given by
    \begin{equation} \label{eq:mingen}
  \min_{\tmmathbf{d}_1 (s), \tmmathbf{d}_2 (s) \in \vv_0^{\bot} (s)}
    \tmmathbf{g} (\tmmathbf{d}_1 (s), \tmmathbf{d}_2 (s)) = \left( \Aop (\vv_0)
  \tmmathbf{n} \cdot \tmmathbf{v}_0 \right)^2 + \left( \Aop \left( \vv_0
  \right) \tmmathbf{b} \cdot \tmmathbf{v}_0 \right)^2,
     \end{equation}
and is reached when
  \begin{align}
    \tmmathbf{d}_1 (\tmmathbf{v}_0) & \assign \tmmathbf{v}_0 \times \left(
     \tmmathbf{v}_0 \times \Aop (\vv_0)
    \tmmathbf{n} \right),  \label{eq:exprdvo1}\\
    \tmmathbf{d}_2 (\tmmathbf{v}_0) & =  \tmmathbf{v}_0 \times \left(
     \tmmathbf{v}_0 \times \Aop (\vv_0)
    \tmmathbf{b} \right) .  \label{eq:exprdvo2}
  \end{align}
  Indeed, one can observe that
  \begin{align}
    \tmmathbf{g} (\tmmathbf{d}_1, \tmmathbf{d}_2) & =  \left| \tmmathbf{d}_1
    -\tmmathbf{v}_0 \times \left(  \tmmathbf{v}_0 \times
    \Aop (\vv_0) \tmmathbf{n} \right) + \left( \Aop \left( \vv_0
    \right) \tmmathbf{n} \cdot \tmmathbf{v}_0 \right) \tmmathbf{v}_0 \right|^2
    \nonumber\\
    &   \qquad \qquad + \left| \tmmathbf{d}_2 -\tmmathbf{v}_0 \times \left(
     \tmmathbf{v}_0 \times \Aop (\vv_0)
    \tmmathbf{b} \right) + \left( \Aop (\vv_0) \tmmathbf{b} \cdot
    \tmmathbf{v}_0 \right) \tmmathbf{v}_0 \right|^2 \nonumber\\
    & =  \left| \tmmathbf{d}_1 -\tmmathbf{v}_0 \times \left(
     \tmmathbf{v}_0 \times \Aop (\vv_0)
    \tmmathbf{n} \right) \right|^2 + \left| \tmmathbf{d}_2 -\tmmathbf{v}_0
    \times \left(  \tmmathbf{v}_0 \times \Aop \left(
    \vv_0 \right) \tmmathbf{b} \right) \right|^2 \nonumber\\
    &   \qquad \qquad + \left( \Aop (\vv_0) \tmmathbf{n} \cdot
    \tmmathbf{v}_0 \right)^2 + \left( \Aop (\vv_0) \tmmathbf{b}
    \cdot \tmmathbf{v}_0 \right)^2 . \nonumber
  \end{align}
  Summarizing, we obtained that
  \begin{align}
    \liminf_{\varepsilon \rightarrow 0}
    \mathcal{E}_{\mathcal{C}}^{\varepsilon} (\tmmathbf{v}_{\varepsilon}) &
    \geqslant  \frac{1}{2} \int_I \left| \partial_s \vv (s) + \Aop \left( \vv
    (s) \right) \tmmathbf{t} (s) \right|^2 \mathd s \nonumber\\
    &   \qquad + \frac{1}{2} \int_I ( \Aop^{\T} (\vv (s))
    \vv (s) \cdot \tmmathbf{n} (s))^2 + ( \Aop^{\T} \left( \vv (s)
    \right) \vv (s) \cdot \tmmathbf{b} (s) )^2 \mathd s, 
  \end{align}
  and this is nothing but the right-hand side of \eqref{eq:gleF}.
  
  \smallskip
  
 {\noindent}{\tmname{($\Gamma$-limsup)}}  
  A direct computation shows that for any $\vv_0\in H^1( \mathcal{C}, \jM)$ such that 
  ${\nabla_{\zv}}  \vv_0 = 0$, the family
  \begin{equation}\label{eq:constrrecseq}
    \tmmathbf{v}_{\varepsilon} \assign \frac{\tmmathbf{v}_0 + \varepsilon
    \tmmathbf{d} (\tmmathbf{v}_0)}{| \tmmathbf{v}_0 + \varepsilon \tmmathbf{d}
    (\tmmathbf{v}_0) |}, \quad \tmmathbf{d} (\tmmathbf{v}_0) \assign z_1
    \tmmathbf{d}_1 (\tmmathbf{v}_0) + z_2 \tmmathbf{d}_2 (\tmmathbf{v}_0),
  \end{equation}
  with  $\tmmathbf{d}_1 (\tmmathbf{v}_0),
  \tmmathbf{d}_2 (\tmmathbf{v}_0)$ given by \eqref{eq:exprdvo1}-\eqref{eq:exprdvo2}, is a recovery sequence, i.e.,
  \begin{equation} \label{eq:recseqtoprove}
  \tmmathbf{v}_{\varepsilon} \rightarrow \tmmathbf{v}_0
  \text{ strongly in } H^1( \mathcal{C}, \jM) \quad \text{ and}
  \quad \mathcal{E}_{\mathcal{C}}^{\varepsilon} (\tmmathbf{v}_{\varepsilon})
  \rightarrow \mathcal{E}_{\mathcal{C}} (\vv_0).
  \end{equation}
  Indeed, since $\tmmathbf{d}_1 (\tmmathbf{v}_0)$, $\tmmathbf{d}_2 (\tmmathbf{v}_0)$ depend only on the $s$-variable and $|\tmmathbf{v}_0)|=1$, we have
  \begin{align}
    D\tmmathbf{v}_{\varepsilon} & = \frac{1}{| \tmmathbf{v}_0 + \varepsilon
    \tmmathbf{d} (\tmmathbf{v}_0) |} \left( I + \frac{\tmmathbf{v}_0 +
    \varepsilon \tmmathbf{d} (\tmmathbf{v}_0)}{| \tmmathbf{v}_0 + \varepsilon
    \tmmathbf{d} (\tmmathbf{v}_0) |} \otimes \frac{\tmmathbf{v}_0 +
    \varepsilon \tmmathbf{d} (\tmmathbf{v}_0)}{| \tmmathbf{v}_0 + \varepsilon
    \tmmathbf{d} (\tmmathbf{v}_0) |} \right) D (\tmmathbf{v}_0 + \varepsilon
    \tmmathbf{d} (\tmmathbf{v}_0)) \nonumber\\
    & \approx  (I +\tmmathbf{v}_0 \otimes \tmmathbf{v}_0) (\partial_s  |
    \partial_1  | \partial_2) (\tmmathbf{v}_0 + \varepsilon \tmmathbf{d}
    (\tmmathbf{v}_0)) \nonumber\\
    & =  (I +\tmmathbf{v}_0 \otimes \tmmathbf{v}_0) (\partial_s
    \tmmathbf{v}_0 + \varepsilon \partial_s \tmmathbf{d} (\tmmathbf{v}_0) |
    \varepsilon \tmmathbf{d}_1 (\tmmathbf{v}_0) | \varepsilon \tmmathbf{d}_2
    (\tmmathbf{v}_0)) \nonumber
  \end{align}
  where we used the approximation symbol to denote equality up to negligible terms when $\varepsilon\to 0$. 
  Given that $0 = \partial_s \tmmathbf{v}_0 \cdot \tmmathbf{v}_0
  =\tmmathbf{d}_1 (\tmmathbf{v}_0) \cdot \tmmathbf{v}_0 =\tmmathbf{d}_2
  (\tmmathbf{v}_0) \cdot \tmmathbf{v}_0$, we infer the following limit relations
  \begin{align}
    \partial_s \tmmathbf{v}_{\varepsilon} & \rightarrow  (I +\tmmathbf{v}_0
    \otimes \tmmathbf{v}_0) \partial_s \tmmathbf{v}_0 = \partial_s
    \tmmathbf{v}_0, \\
    \varepsilon^{- 1} \partial_1 \tmmathbf{v}_{\varepsilon} & \rightarrow (I
    +\tmmathbf{v}_0 \otimes \tmmathbf{v}_0) \tmmathbf{d}_1 (\tmmathbf{v}_0)
    =\tmmathbf{d}_1 (\tmmathbf{v}_0), \\
    \varepsilon^{- 1} \partial_2 \tmmathbf{v}_{\varepsilon} & \rightarrow  (I
    +\tmmathbf{v}_0 \otimes \tmmathbf{v}_0) \tmmathbf{d}_2 (\tmmathbf{v}_0)
    =\tmmathbf{d}_2 (\tmmathbf{v}_0) ,
  \end{align}
  from which \eqref{eq:recseqtoprove} follows.
  
\smallskip
  
 {\noindent}{\tmname{(Strong convergence of minimizers)}} Let $(\vv_{\varepsilon})_{\varepsilon} \subset H^1 (\mathcal{C},
\mathbb{S}^2)$ be a minimizing family for
$(\mathcal{E}_{\mathcal{C}}^{\varepsilon})_{\varepsilon}$. By compactness (cf.
\eqref{eqs:compact1}--\eqref{eqs:compact2}), there exists $\vv_0 \in H^1
(\mathcal{C}, \mathbb{S}^2)$, depending only the $s$-variable, such that
$\vv_{\varepsilon} \rightharpoonup \vv_0$ weakly in $H^1 (\mathcal{C},
\mathbb{S}^2)$ and $\partial_s \vv_{\varepsilon} \rightarrow 0$ strongly in
$L^2 (\mathcal{C}, \mathbb{S}^2)$. By the $\Gamma$-liminf inequality, the
strong convergence of $\vv_{\varepsilon} \rightarrow \vv_0$ in $L^2
(\mathcal{C}, \mathbb{S}^2)$, and the minimality of $\vv_{\varepsilon}$, we
get that if $\vv_{\varepsilon}^{\star}$ is the recovery sequence built from
$\vv_0$ as in \eqref{eq:constrrecseq}, then
\begin{equation}
  \mathcal{E}_{\mathcal{C}} (\vv_0) \leqslant
  \liminf_{\varepsilon \rightarrow 0} \mathcal{E}_{\mathcal{C}}^{\varepsilon}
  (\vv_{\varepsilon}) \leqslant \limsup_{\varepsilon \rightarrow
  0} \mathcal{E}_{\mathcal{C}}^{\varepsilon} (\vv_{\varepsilon})
  \leqslant \lim_{\varepsilon \rightarrow 0}
  \mathcal{E}_{\mathcal{C}}^{\varepsilon} (\vv_{\varepsilon}^{\star}) \eqs \mathcal{E}_{\mathcal{C}} (\vv_0) .
  \label{eq:limrelforstrongconv}
\end{equation}
Therefore, if $(\vv_{\varepsilon})_{\varepsilon}$ is a minimizing
family for $(\mathcal{E}_{\mathcal{C}}^{\varepsilon})_{\varepsilon}$, then
$\vv_{\varepsilon} \rightharpoonup \vv_0$ weakly in $H^1 (\mathcal{C},
\mathbb{S}^2)$ for some $\vv_0$ which depends only on the $s$-variable and,
moreover, by \eqref{eq:limrelforstrongconv}, $\lim_{\varepsilon \rightarrow 0}
\mathcal{E}_{\mathcal{C}}^{\varepsilon} (\vv_{\varepsilon})
=\mathcal{E}_{\mathcal{C}} (\vv_0)$. By the compactness relations
\eqref{eqs:compact1}--\eqref{eqs:compact2} and the lower semicontinuity of the
norm, we get that
\begin{align}
  \mathcal{E}_{\mathcal{C}} (\vv_0) = \lim_{\varepsilon
  \rightarrow 0} \mathcal{E}_{\mathcal{C}}^{\varepsilon}(\vv_{\varepsilon}) & \geqslant  \limsup_{\varepsilon \rightarrow 0}
  \frac{1}{2 |Q|}  \int_{I \times Q} \left| \partial_s \vv_{\varepsilon} +
  \Aop (\vv_{\varepsilon}) \tmmathbf{t} \right|^2 \mathd \zv
  \mathd s \nonumber\\
  &  \qquad \qquad + \frac{1}{2 |Q|}  \int_{I \times Q}
  \tmmathbf{g} (\tmmathbf{d}_1, \tmmathbf{d}_2) \mathd \zv \mathd s
  \nonumber\\
  & \geqslant \limsup_{\varepsilon \rightarrow 0}  \frac{1}{2 |Q|}  \int_{I
  \times Q} \left| \partial_s \vv_{\varepsilon} + \Aop \left(
  \vv_{\varepsilon} \right) \tmmathbf{t} \right|^2 \mathd \zv \mathd s
  \nonumber\\
  &   \qquad \qquad  + \frac{1}{2}  \int_I \left( \Aop^{\T}
  (\vv_0) \vv_0 \cdot \tmmathbf{n} \right)^2 + \left( \Aop^{\T}
  (\vv_0) \vv_0 \cdot \tmmathbf{b} \right)^2 \mathd s, 
  \label{eqs:templimsupstrong}
\end{align}
where, for the last inequality, we used the fact that the function
$\tmmathbf{g}$ defined by \eqref{eq:exprg} has minimal energy $\left( \Aop
(\vv_0) \tmmathbf{n} \cdot \tmmathbf{v}_0 \right)^2 + \left( \Aop
(\vv_0) \tmmathbf{b} \cdot \tmmathbf{v}_0 \right)^2$.
Substituting in \eqref{eqs:templimsupstrong} the equivalent expression of
$\mathcal{E}_{\mathcal{C}} (\vv_0)$ given in
Remark~\ref{rmk:altdescgammlim}, we get that
\begin{equation}
  \limsup_{\varepsilon \rightarrow 0}  \frac{1}{2 |Q|}  \int_{I \times Q}
  \left| \partial_s \vv_{\varepsilon} + \Aop (\vv_{\varepsilon})
  \tmmathbf{t} \right|^2 \mathd \zv \mathd s \leqslant \frac{1}{2}  \int_I
  \left| \partial_s \vv_0 + \Aop (\vv_0) \tmmathbf{t} \right|^2
  \mathd s.
\end{equation}
Hence, from the lower semicontinuity of the norm, it follows that
\begin{equation}
  \lim_{\varepsilon \rightarrow 0}  \frac{1}{2 |Q|}  \int_{I \times Q} \left|
  \partial_s \vv_{\varepsilon} + \Aop (\vv_{\varepsilon})
  \tmmathbf{t} \right|^2 \mathd \zv \mathd s = \frac{1}{2}  \int_I \left|
  \partial_s \vv_0 + \Aop (\vv_0) \tmmathbf{t} \right|^2 \mathd
  s.
\end{equation}
The previous relation implies that $\left\| \partial_s
\vv_{\varepsilon} \right\|_{L^2 (\mathcal{C}, \mathbb{S}^2)} \rightarrow
\left\| \partial_s \vv_0 \right\|_{L^2 (\mathcal{C}, \mathbb{S}^2)}$. Overall,
also considered that $\left\| \partial_1 \vv_{\varepsilon} \right\|_{L^2
(\mathcal{C}, \mathbb{S}^2)} \rightarrow 0$ and $\left\| \partial_2
\vv_{\varepsilon} \right\|_{L^2 (\mathcal{C}, \mathbb{S}^2)} \rightarrow 0$,
we deduce the convergence of the norms $\left\| \vv_{\varepsilon}
\right\|_{H^1 (\mathcal{C}, \mathbb{S}^2)} \rightarrow \left\| \vv_0
\right\|_{H^1 (\mathcal{C}, \mathbb{S}^2)}$ which together with the weak
convergence $\vv_{\varepsilon} \rightharpoonup \vv_0$ in $H^1 (\mathcal{C},
\mathbb{S}^2)$, assures strong convergence in $H^1 (\mathcal{C},
\mathbb{S}^2)$. This concludes the proof.
\end{proof}

\section{Applications to Micromagnetics}
In this section, we illustrate how our analysis of curved thin wires can
explain different scenarios that may occur when ferromagnetic crystals lack
inversion symmetry and how the so-called Dzyaloshinskii-Moriya interaction
(DMI) can distort the otherwise uniform ferromagnetic state. Before delving
into specifics, we describe the physical framework and introduce relevant
notation.

In the variational theory of micromagnetism ({cf.}~{\cite{BrownB1963,hubert2008magnetic,LandauA1935}}), the order parameter that
characterizes the orientations of a rigid ferromagnetic body within a region
$\Omega \subseteq \mathbb{R}^3$ is the magnetization $\tmmathbf{M}$, a vector
field whose magnitude $M_s \assign | \tmmathbf{M} |$, called the saturation
magnetization, can be considered constant in $\Omega$ at temperatures well
below the Curie temperature. Hence, one can represent the magnetization as
$\tmmathbf{M} \assign M_s \vv$, where $\vv : \Om \rightarrow \mathbb{S}^2$ is
a map taking values into the unit sphere of $\mathbb{R}^3$. Even though the magnitude of $\vv$ remains constant in space, its direction is free to vary and the observable magnetization states can be expressed as local minimizers
of the micromagnetic energy functional. In single-crystal ferromagnets
({cf.}~{\cite{alouges2015homogenization,Davoli_2020}}),
and upon appropriate non-dimensionalization, this functional reads as
\begin{equation}
  \mathcal{E}_{\Omega} (\tmmathbf{v}) +\mathcal{W}_{\Omega} (\tmmathbf{v})
  \assign \frac{1}{2}  \int_{\Omega} | \nabla \tmmathbf{v}|^2 \mathd x +
  \frac{1}{2}  \int_{\Omega} | \nabla u |^2 \mathd x \label{eq:GLunorm}
\end{equation}
with $\tmmathbf{v} \in H^1 (\Omega, \mathbb{S}^2)$, and where $\tmmathbf{v}
\chi_{\Omega}$ denotes the extension of $\tmmathbf{v}$ by zero to the whole
space outside $\Omega$.
In the above expression, the exchange energy $\mathcal{E}_{\Omega}$ accounts
for {\tmem{symmetric}} exchange interactions, which tend to penalize
nonuniformities in the magnetization orientation. The term
$\mathcal{W}_{\Omega}$ represents the self-energy due to magnetostatics and
describes the energy due to the magnetic scalar potential $u$
produced by $\tmmathbf{m}$. For $\Om$ bounded, the magnetic scalar potential
can be defined as the unique solution in $H^1 (\mathbb{R}^3)$ of the Poisson
equation 
\begin{equation}
\Delta u = \nabla \cdot \left( \vv \chi_{\Om_{\varepsilon}}
\right)
\end{equation} 
We refer the reader to \cite{doi:10.1137/19M1261365} for further mathematical details.


This paper goes beyond the classical {\tmem{symmetric}} exchange interactions
by examining the potential lack of centrosymmetry in the crystal lattice
structure of the ferromagnet. In this paper, other than the classical
{\tmem{symmetric}} exchange interactions, we consider the possible lack of
centrosymmetry in the crystal lattice structure of the ferromagnet. In
addition to the energy density expressed in equation \eqref{eq:GLunorm}, we
consider possible {\tmem{antisymmetric}} exchange interactions by including
contributions from relevant Lifshitz invariants of the chirality tensor
$\nabla \vv \times \vv$. Of particular significance is the energy contribution
from {\tmem{bulk}} Dzyaloshinskii-Moriya interaction (DMI), with energy
density $\vv \cdot \left( \nabla \times \vv \right)$ corresponding to the
trace of the chirality tensor. Specifically, for every $\vv \in H^1 (\Om,
\mathbb{S}^2)$, we define the bulk DMI energy as
\begin{equation}
  \mathcal{H}_{\Omega} (\vv) \assign \kappa \int_{\Omega} \vv \cdot \left(
  \nabla \times \vv \right) \hspace{0.17em} \mathd x.
\end{equation}
The normalized constant $\kappa \in \mathbb{R}$ is the so-called DMI constant,
and its sign determines the chirality of the system.

The full micromagnetic
energy functional we are interested in is then $\mathcal{E}_{\Omega}
+\mathcal{H}_{\Omega} +\mathcal{W}_{\Omega}$.
In particular, with the notation introduced in Section~\ref{sec:sec2}, the
micromagnetic energy associated with an $\varepsilon$-tube $\Om_{\varepsilon}$
along a curve $\gamma \of I \mapsto \mathbb{R}^3$ with cross-sectional shape
$Q$ is given by
\begin{equation}
\mathcal{G}_{\varepsilon}( \vv ) \assign \frac{1}{\varepsilon^2
   | Q |} \left( \frac{1}{2} \int_{\Om_{\varepsilon}}  \left| D \vv
    \right|^2 \mathd x + \kappa \int_{\Om_{\varepsilon}}  \vv\cdot \nabla \times \vv  
   \, \mathd x + \frac{1}{2} \int_{\mathbb{R}^3} | \nabla u |^2 \mathd x
   \right), 
\end{equation}
We note that in the energy above, we take a form of so-called {\it interior DMI} $\nabla \times \vv \cdot \vv$. To put it in the form \eqref{eq:mainenfunc}, we define a matrix $K$ as follows
\begin{equation}
K( \vv ) = \kappa \left(\begin{array}{ccc}
      0 & - v_3 & v_2\\
     v_3 & 0 & - v_1\\
     - v_2 & v_1 & 0
   \end{array}\right).\notag
\end{equation}
It is clear that
\[ \frac{1}{2} \left| D \vv + K ( \vv ) \right|^2 = \frac{1}{2}
   \left| D \vv \right|^2 + D \vv : K ( \vv) + \frac{1}{2} \left|
   K( \vv ) \right|^2, \]
and a simple computation shows that $ D \vv : K ( \vv ) = \kappa \vv\cdot \nabla \times \vv$, and
$| K ( \vv) |^2 = 2 \kappa^2$.
Based on these relations, we can rewrite our energy as follows
\[ \mathcal{G}_{\varepsilon}( \vv ) \assign \frac{1}{\varepsilon^2
   | Q |} \left(\frac{1}{2} \int_{\Om_{\varepsilon}}  \left| D \vv
    + K( \vv) \right|^2 \mathd x + \frac{1}{2}
   \int_{\mathbb{R}^3} | \nabla u |^2 \mathd x \right) - \kappa^2 \frac{|
   \Omega_{\varepsilon} |}{\varepsilon^2 | Q |} . \]
Without loss of generality, we assume $| Q | = 1$ and note that the last constant term in the energy can be dropped as it does not affect critical points or $\Gamma$-convergence. Then, using $\Gamma$-convergence results from \cite{Slastikov2012} (see also \cite{Harutyunyan2016}) and \eqref{eq:gleF}, in the limit as $\varepsilon \to 0$ we obtain the
following limiting energy
\begin{equation} \label{eq:fulllimfunc}
\mathcal{E}_{\mathcal{C}} ( \vv ) = \frac{1}{2} \int_I \left| \partial_s \vv (s) + \Aop ( \vv (s)) \tmmathbf{t} (s) \right|^2 \mathd s + \frac{1}{2} \int_I M \vv (s) \cdot \vv (s) \mathd s ,
\end{equation}
where
\begin{equation} \label{eq:fulllimfunc2} 
M = - \frac{1}{2 \pi} \int_{\partial Q} \int_{\partial Q} \tmmathbf{n}(\xi) \otimes \tmmathbf{n} (\eta) \ln |\xi-\eta| \mathd \xi \mathd \eta . 
\end{equation}
We now want to investigate several special cases of the limiting energy and find its minimizers in the form of domain walls. 

\subsection{Straight wire with DMI}
In the case of an infinite straight wire with a circular cross-section, the interval $I$ in \eqref{eq:fulllimfunc} is not compact, and one cannot directly apply Theorem~\ref{thm:main}. However, it is still possible to prove the result for an infinite straight wire with some standard modifications to account for translation invariance. It is also possible to model this scenario by considering it as the limit of a sequence of energy functionals defined on straight wires of finite lengths. We have $\tmmathbf{t}=\tmmathbf{e}_1$ and, up to constant terms, the limiting energy can be rewritten as
\begin{equation}\label{eq:strwireen}
    \mathcal{E}_{\mathcal{C}}( \vv) =
    \frac{1}{2} \int_{\RR} \left| \partial_x \vv  + \kappa(0,v_3,-v_2)\right|^2 \mathd x + \frac{1}{4\pi} \int_{\RR}
     (1-v_1^2) \mathd x \, .
\end{equation} 
We note that the global minimum of the energy is achieved for $\vv=\pm \tmmathbf{e}_1$.

We are interested in non-trivial local minimizers corresponding to domain walls. Therefore, we can impose boundary conditions $\vv \to \pm \tmmathbf{e}_1$ as $x\to \pm \infty$ and look for a global minimizer of $\mathcal{E}_{\mathcal{C}}$ with these boundary conditions.  For that, we decompose the first term as 
\begin{equation}
    \left| \partial_x \vv + \kappa(0,v_3,-v_2)\right|^2= |v_1'|^2+ |v_2'+\kappa v_3|^2+|v_3'-\kappa v_2|^2.
\end{equation}
and represent the magnetization $\vv$ in spherical coordinates
\[
\vv(x)=(\cos\theta(x), \sin\theta(x)\cos\phi(x),\sin\theta(x)\sin\phi(x))\, .
\]
Plugging it into the energy we obtain
\begin{equation}
\mathcal{E_C}(\theta,\phi)=\frac12 \int_\RR |\theta'(x)|^2 + \sin^2\theta(x) |\phi'(x)-\kappa|^2 \, \mathd x + \frac{1}{4\pi} \int_\RR \sin^2\theta(x) \, \mathd x.
\end{equation}
The boundary conditions translate into $\theta(x) \to 0$ as $x\to \infty$,  $\theta(x) \to k\pi$ as $x\to -\infty$, where $k\in \NN$ is odd.
We minimize this energy in $\phi$ to obtain $\phi(x)=\kappa x$. The subsequent minimization in $\theta$ gives us $k=1$ and a standard domain wall profile
\begin{align}
    \theta_0(x)=2\arctan\left(e^{-\frac{1}{\sqrt{4\pi}} x} \right).
\end{align}
The minimizing profile is
\begin{align}
\vv_0(x)=(\cos\theta_0(x), \sin\theta_0(x)\cos(\kappa x), \sin\theta_0(x)\sin(\kappa x)).
\end{align}
It is clear that $\vv_0$ is a global minimizer of the energy \eqref{eq:strwireen} under domain wall boundary conditions $\vv \to \pm \tmmathbf{e}_1$ as $x\to \pm \infty$. Therefore, it follows that it is a local minimizer of $\mathcal{E}_{\mathcal{C}}$ with respect to $H^1$ perturbations.

\begin{figure*}[t]
	\centering
	\includegraphics[width=15.4cm]{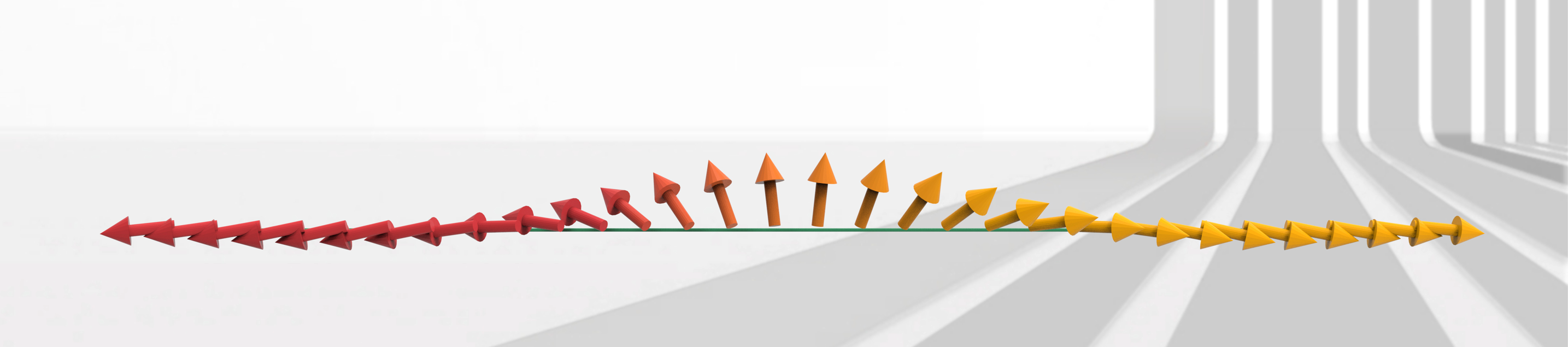}
	\caption{\small
Sketch of the domain wall in the straight wire. 
The value of the dimensionless parameter characterizing the strength of the DMI is chosen to be $\kappa=0.36$, which is close to that for the B20-type FeGe~\cite{Zheng2021}.
}
	\label{fig_dw}
\end{figure*}

We note that in the intuitive 1D energy for a nanowire sometimes studied in physics literature for soft materials \cite{boulle13}
\begin{align}
  \frac{1}{2} \int_{\RR} \left| \partial_x \vv\right|^2 + \kappa(v_2' v_3-v_3' v_2)  +Q (1-v_1^2)  \, \mathd x
\end{align}
 the DMI term might lead to an uncontrolled profile winding at the tails to decrease the micromagnetic energy.  This phenomenon is always absent in the reduced 3D to 1D model due to the presence of two effective anisotropy terms produced by DMI and stray field energies.
Effective anisotropy here globally suppresses conical modulation~\cite{Borisov1985}, but modulation occurs along the domain wall (see Fig.~\ref{fig_dw}).


\subsection{Rings with DMI} If our original domain has a geometry of a ring with radius $R$ and circular cross-section, we can define $\tmmathbf{t}
(\phi) = (- \sin \phi, \cos \phi, 0)$, $\tmmathbf{n} (\phi) = (\cos \phi, \sin
\phi, 0)$ and hence
\[ \Aop ( \vv (\phi)) \tmmathbf{t} (\phi) = \kappa (- v_3 \cos
   \phi, - v_3 \sin \phi, v_2 \sin \phi + v_1 \cos \phi) = - \kappa v_3
   \tmmathbf{n} (\phi) + \kappa \vv (\phi) \cdot \tmmathbf{n} (\phi)
   \tmmathbf{e}_3 . \]
It is convenient to rescale the domain to be a unit circle, and hence, we can write the energy rescaled by $R$ as
\begin{align}
\mathcal{E}_{\mathcal{C}} ( \vv ) = 
    \displaystyle \frac{1}{2R^2} \int_0^{2 \pi} \left| \partial_{\phi} \vv  \right|^2
     \mathd \phi + \frac{\kappa}{R} \int_0^{2 \pi} \left( - \partial_{\phi} \vv
     \cdot \tmmathbf{n}  v_3 + \vv  \cdot \tmmathbf{n} 
     \partial_{\phi} v_3 \right) \mathd \phi \\
     \hspace{10em} + \displaystyle \frac{\kappa^2}{2} \int_0^{2 \pi} 1 - \left( \vv
      \cdot \tmmathbf{t}  \right)^2 \mathd \phi + \frac{1}{4\pi}
     \int_0^{2 \pi}  1 - \left( \vv
      \cdot \tmmathbf{t}  \right)^2 \mathd \phi .
\end{align}
When DMI can be neglected, an extensive examination of energy minimizers and the resulting symmetry-breaking phenomena is carried out in \cite{DifFioSla2023} (see also \cite{CarbouRing2022} for findings on dynamic stability).

\begin{figure*}[t]
	\centering
	\includegraphics[width=15.4cm]{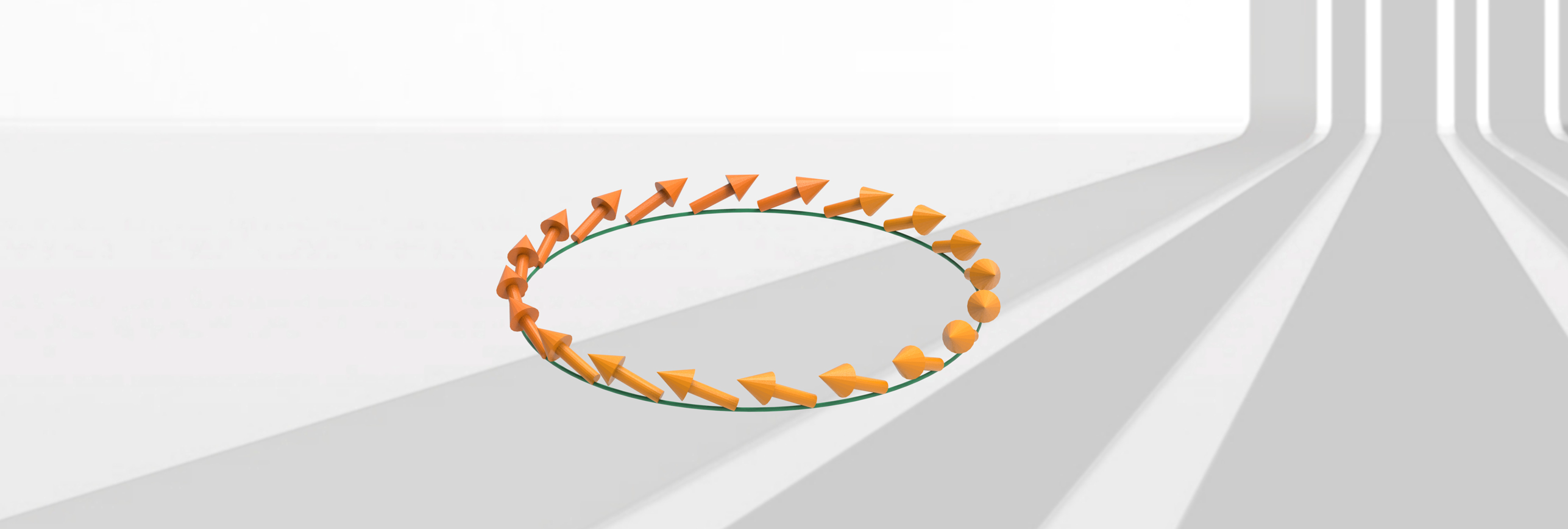}
	\caption{\small
Sketch of the minimizer in the circular wire. 
Spin configuration is drawn according to equation~(\ref{sol:trivialwithaniso}) with an arbitrarily chosen value of $\gamma=2.3$.
}\label{fig_ring}
\end{figure*}

Observe that the DMI term produces an additional anisotropy along the wire. It is convenient to introduce local ring coordinates and decompose the magnetization $\vv$ along the Frenet frame $\vv=v_t \tmmathbf t + v_n \tmmathbf n+ v_3 {\tmmathbf{b}}$. We now have $\partial_\phi \vv = (v_t' + v_n)\tmmathbf t + (v_n' - v_t)\tmmathbf n + v_3' {\tmmathbf{b}}$ and, therefore, 
it is possible to rewrite the reduced 1D energy as
\begin{align}
   \mathcal{E}_{\mathcal{C}} ( \vv )  & =    \frac{1}{2R^2} \int_0^{2 \pi} \left|v_t' + v_n\right|^2 + \left|v_n'-v_t -R\kappa v_3 \right|^2 +\left|v_3'+R\kappa v_n \right|^2
     \mathd \phi  \notag\\
     & \hspace{20em}  +  \frac{1}{4\pi}
     \int_0^{2 \pi} (1-v_t^2)\mathd \phi 
\end{align}
Let us first try to understand the influence of the DMI in the model with only exchange and DMI present 
\[ \mathcal{E}^0_{\mathcal{C}} ( \vv ) = 
     \displaystyle \frac{1}{2R^2} \int_0^{2 \pi} \left|v_t' + v_n\right|^2 + \left|v_n'-v_t -R\kappa v_3 \right|^2 +\left|v_3'+R\kappa v_n \right|^2
     \mathd \phi .
     \]
It is clear that exchange energy wants magnetization to be a constant vector. However, the presence of DMI breaks the symmetry of the problem and leads to novel types of minimizers. Indeed, we can solve the system
\begin{align}
    v_t'+v_n=0,\qquad v_n'-v_t-R\kappa v_3 =0, \qquad v_3'+R\kappa v_n=0, \qquad |\vv|=1,
\end{align}
and obtain the following set of solutions:
\begin{align}
    v_t&=\frac{A}{\sqrt{1+R^2\kappa^2}} \cos(\sqrt{1+R^2\kappa^2}\phi {+ \phi_0}) -R\kappa B,\label{eq:exsol01}\\
    v_n&= A\sin(\sqrt{1+R^2\kappa^2}\phi {+ \phi_0}), \label{eq:exsol02} \\ 
    v_3 &= \frac{A R\kappa}{\sqrt{1+R^2\kappa^2}} \cos(\sqrt{1+R^2\kappa^2}\phi {+ \phi_0}) +B, \label{eq:exsol}
\end{align}
where
\begin{align}\label{con:ab}
A^2+B^2(1+R^2\kappa^2)=1.
\end{align}
We also need to make sure that periodicity conditions on $\vv$ are satisfied and therefore, provided $A \neq 0$, we have the following condition
\begin{align}\label{con:quant}
(1+R^2\kappa^2)=n^2,
\end{align}
where $n\in \NN$. Note that all these solutions yield zero micromagnetic energy (with absent stray field), and hence, depending on the relation between $R$ and $\kappa$, we might have several minimizers. Two minimizers always exists by taking $A=0$ and $B=\pm\frac{1}{\sqrt{1+R^2\kappa^2}}$ (see figure~\ref{fig_ring}).

It is interesting to observe that the introduction of DMI breaks the degeneracy of exchange interaction, and instead of any constant solution being a minimizer of symmetric exchange energy we obtain only two global minimizers of energy $\mathcal{E}^0_{\mathcal{C}}$ with anti-symmetric exchange 
\begin{align}
    \vv(\phi)=\mp \frac{R \kappa}{\sqrt{1+R^2\kappa^2}} \tmmathbf t \pm \frac{1}{\sqrt{1+R^2\kappa^2}} {\tmmathbf{b}}.
\end{align}
Moreover, when condition \eqref{con:quant} is satisfied we have a continuum of solutions (cf.~\eqref{eq:exsol01}--\eqref{eq:exsol}) with $A$ and $B$ satisfying \eqref{con:ab}.

An alternative approach to solving the problem sheds light on non-trivial periodic solutions beyond condition~(\ref{con:quant}). 
Let's define an auxiliary unit vector field as follows:
\begin{align}
\tilde{\vv} = 
\left(\begin{array}{ccc}
     1 & 0 & 0\\
     0 & \cos(\alpha) & \sin(\alpha)\\
     0 & -\sin(\alpha) & \cos(\alpha)
   \end{array}\right) \cdot
\left(\begin{array}{ccc}
     \cos(\phi) & \sin(\phi) & 0\\
     -\sin(\phi) & \cos(\phi) & 0\\
     0 & 0 & 1
   \end{array}\right) \cdot
   \vv,
\end{align}
where $\alpha = \arctan(R\kappa)$. 
The expression for the energy depending on field $\tilde{\vv}$ has the following form:
\[  \begin{array}{ll}
    \displaystyle \mathcal{E}^0_{\mathcal{C}} ( \tilde{\vv} ) = \frac{1}{2R^2} \int_0^{2 \pi} \left| \partial_{\phi} \tilde{\vv}  \right|^2
     \mathd \phi + 
     \frac{\sqrt{1+R^2\kappa^2}}{R^2} \int_0^{2 \pi} 
     \left(  
     \tilde{v}_1 \partial_\phi \tilde{v}_2 -  \tilde{v}_2 \partial_\phi \tilde{v}_1
     \right) \mathd \phi &  \\
     \hspace{10em} \hspace{10em} + \displaystyle \frac{1+R^2\kappa^2}{2 R^2} \int_0^{2 \pi} ( 1 - \tilde{v}_3^2)  \mathd \phi \, .
   \end{array}
   \]
The resulting Hamiltonian clearly shows that the system is frustrated between two states: a homogeneous state polarized along the easy axis or a N\'{e}el spiral in a perpendicular plane. 
We note that it is possible to represent $\tilde \vv$ in spherical coordinates as follows
$$
\tilde \vv = (\sin\theta(\phi)\cos\psi(\phi),\sin\theta(\phi)\sin\psi(\phi),\cos\theta(\phi))
$$
and rewrite the above energy as
\begin{align}
\mathcal{E}^0_{\mathcal C}(\theta,\phi)=\frac{1}{2R^2} \int_0^{2\pi} |\theta'|^2 + \sin^2\theta |\psi'-\sqrt{1+R^2\kappa^2}|^2 \, \mathd \phi  .
\end{align}
We already have information about minimizers of the energy, but using this representation, we can also find special solutions to the Euler-Lagrange equations 
\begin{align}
-\theta''+\sin\theta \cos\theta |\psi'-\sqrt{1+R^2\kappa^2}|^2 &=0\\
-(\sin^2\theta \psi')' + \sqrt{1+R^2\kappa^2}(\sin^2 \theta)'&=0
\end{align}
Here we have to supplement these with boundary conditions $\theta(0)=\theta(2\pi)+2k\pi$, $\psi(0)=\psi(2\pi)+2m\pi$, $\theta'(0)=\theta'(2\pi)$, $\psi'(0)=\psi'(2\pi)$.
We can take $\theta=\frac{\pi}{2}$ and $\psi=n \phi + \phi_0$, $n \in \NN$ to deduce
\begin{align} 
\tilde{\vv} & =  (\sin(n \phi + \phi_0), \cos(n \phi + \phi_0), 0). \label{sol:II}
\end{align}
Solution~(\ref{sol:II}) is almost always a saddle point, except when condition~(\ref{con:quant}) is satisfied. 
With this exceptional condition, the unstable mode disappears, and the solution~(\ref{sol:II}) coincides with the minimizer.

If we introduce the magnetostatic energy contribution and study energy $\mathcal{E_C}$, in general, the minimization problem cannot be solved explicitly anymore, and numerical simulations are needed to investigate the minimizing configurations. However, we can still find an explicit solution  
\begin{align}
    \vv(\phi) &=\mp \frac{\gamma}{\sqrt{1+\gamma^2}} \tmmathbf t \pm \frac{1}{\sqrt{1+\gamma^2}} \tmmathbf{b}, \label{sol:trivialwithaniso} \\ 
    \gamma &= \frac{\sqrt{4R^2\kappa^2 + (R^2\kappa^2 + \tfrac{1}{2\pi} R^2 - 1)^2} + R^2\kappa^2 + \tfrac{1}{2\pi} R^2 - 1}{2R\kappa}.
\end{align}
It corresponds to taking $\theta$ and $\psi$ as constants and solving an algebraic problem coming from Euler-Lagrange equations for $\theta$ and $\psi$. Our additional numerical simulations indicate that \eqref{sol:trivialwithaniso} is a good candidate for minimizer in a wide range of values $R$ and $\kappa$. Still, the rigorous clarification of this question is beyond the scope of this work.

\begin{remark} We note that our results also apply to some interactions beyond exchange and DMI. Below, we give one example of such a situation --- the Ado interaction --- highly nonlinear interaction possible in several systems with broken inversion symmetry~\cite{Rybakov2021, Ado2021}.
According to~\cite{Rybakov2021, Ado2021}, for the crystals with the point group symmetry $T_\text{d}$ we have the following simplified energy
\begin{equation}
\mathcal{G}_{\varepsilon} ( \vv ) \assign \frac{1}{\varepsilon^2
   | Q |} \left( \int_{\Om_{\varepsilon}} \frac{1}{2} \left| D \vv
    \right|^2 + \beta (v_1
   v_2  v_3)\nabla \cdot \vv   \, \mathd x \right).
\end{equation}
Given that $K ( \vv ) \tmmathbf{e}_1 \cdot \partial_x \vv + K ( \vv
   ) \tmmathbf{e}_2 \cdot \partial_y \vv + K ( \vv )
   \tmmathbf{e}_3 \cdot \partial_z \vv = \beta (v_1
   v_2  v_3)\nabla \cdot \vv$ when $K ( \vv ) = \beta (v_1  v_2  v_3)I$ with $I$ the $3\times 3$ identity matrix, using Theorem~\ref{thm:main} we
can obtain the limiting energy for Ado interactions
\begin{equation}
\mathcal{E}_{\mathcal{C}} ( \vv ) = 
     \frac{1}{2} \int_I \left| \partial_x \vv  + \beta v_1  v_2  v_3
      \tmmathbf{e}_1 \right|^2 \mathd x - \frac{3 \beta^2}{2} \int_I v_1^2
      v_2^2 v_3^2 \mathd x \, .
\end{equation}

\end{remark}

\section*{Acknowledgments}
\noindent
{\tmname{G.Di~F.}} is a member of Gruppo Nazionale per l'Analisi Matematica, la Probabilità e le loro Applicazioni (GNAMPA) of INdAM. {\tmname{G.Di~F.}} acknowledges partial support from the Italian Ministry of Education and Research through the PRIN2022 project {\emph{Variational Analysis of Complex Systems in Material Science, Physics and Biology}} No.~2022HKBF5C, and from the Austrian Science Fund (FWF) through the project {\emph{Analysis and Modeling of Magnetic Skyrmions}} (grant 10.55776/P34609). {\tmname{G.Di~F.}} also thanks TU Wien and MedUni Wien for their hospitality. \tmname{F.N.R.} acknowledges support from the Swedish Research Council (grant 2023-04899).

\bibliographystyle{siam} 
\bibliography{literature}

\end{document}